\documentclass[12pt]{amsart}
\usepackage{amsmath,amsthm}

\def\pmatrix{\left(\begin{matrix}}
\def\endpmatrix{\end{matrix}\right)}

\def\Sp{\operatorname{Sp}}

\def\R{{\mathbb R}}

\def\Z{{\mathbb Z}}
\def\F{{\mathbb F}_2}
\def\C{{\mathbb C}}

\def\Q{{\mathbb Q}}
\def\Pj{{\mathbb P}}

\def\cal{\mathcal}
\def\Cal{\mathcal}

\def\slope{\operatorname{slope}}
\def\Norm{\operatorname{Norm}}
\def\supp{\operatorname{supp}}

\def\diag{\operatorname{diag}}
\def\Stab{\operatorname{Stab}}
\def\GL{\operatorname{GL}}
\def\SL{\operatorname{SL}}

\def\Sym{\operatorname{Sym}}
\def\Jac{\operatorname{Jac}}

\def\de{\delta}

\def\tr{{\rm tr}}
\def\Tr{{\rm Tr}}
\def\inv{^{-1}}
\def\<{\langle}
\def\>{\rangle}

\def\t{\theta}

\def\e{\varepsilon}

\def\a{\alpha}
\def\b{\beta}

\def\A{{\mathcal A}}

\def\H{{\mathcal H}}

\def\tch#1#2{{\left[\begin{matrix}#1\\ #2\end{matrix}\right]}}
\def\tt#1#2{{\t\tch{#1}{#2}}}
\def\smtwomat#1#2#3#4{{\left(\begin{matrix}#1 & #2\\ #3 & #4\end{matrix}\right)}}
\def\Xg{{\mathcal X}_g}
\def\Pg{{\mathcal P}_g}
\def\Pgone{{\mathcal P}_{g_1}}
\def\Pgtwo{{\mathcal P}_{g_2}}
\def\PR{\Pg(\R)}
\def\PsR{\Pg(\R)^{\text{semi}}}

\theoremstyle{plain}
\newtheorem{thm}{Theorem}
\newtheorem{lm}[thm]{Lemma}
\newtheorem{prop}[thm]{Proposition}
\newtheorem{cor}[thm]{Corollary}

\newtheorem{df}[thm]{Definition}
\newtheorem{rem}[thm]{Remark}
\theoremstyle{definition}
\begin{document}
\title[Modular Forms of weight 8]{ Modular Forms of weight $8$ for $\Gamma_g(1,2)$}

\author[M. Oura]{Manabu Oura}
\address{Department of Mathematics, Faculty of Science, Kochi University, 2-5-1 Akebono-cho, Kochi 780-8520, Japan}
\email{oura@kochi-u.ac.jp}

\author[C. Poor]{Cris Poor}
\address{Department of Mathematics, Fordham University, Bronx, NY 10458-5165, USA}
\email{poor@fordham.edu}

\author[R. Salvati Manni]{Riccardo Salvati Manni}
\address{Dipartimento di Matematica, Universit\`a ``Sapienza'',
Piazzale A. Moro 2, Roma, I 00185, Italy}
\email{salvati@mat.uniroma1.it}

\author[D. Yuen]{David S. Yuen}
\address{Department of Mathematics and Computer Science, Lake Forest College, 555 N. Sheridan Rd., Lake Forest, IL 60045, USA}
\email{yuen@lakeforest.edu}

\date{\today}
\begin{abstract}
We complete the program indicated by the Ansatz of D'Hoker and Phong
in genus ~$4$ by proving the uniqueness of the restriction to Jacobians of
the weight~$8$ Siegel cusp
forms  satisfying the Ansatz.  
We prove $\dim [\Gamma_4(1,2),8]_0=2$ and $\dim [\Gamma_4(1,2),8]=7$.  
In each genus, we classify the linear relations among the self-dual lattices
of rank~{16}.  
We extend the program to genus ~$5$ by constructing the unique linear combination
of theta series that satisfies the Ansatz.  
\end{abstract}
\maketitle

\section{Introduction}
Modular forms of weight 8 with respect to the theta group $\Gamma_g(1,2)$  have  recently been  a useful tool in Physics.
Some are  fundamental in the construction of a chiral superstring
measure, cf. \cite{DHP1}, \cite{DHPa}, \cite{DHPb}, \cite{DHPc}  and  \cite{GR}, for small genera.  
Let $\H_g$ be the Siegel upper half space and let the {\it Jacobian locus\/},
$\Jac_g \subseteq \H_g$,
be the set of period matrices of compact Riemann surfaces.  
For $g \le 3$, $\Jac_g$ is dense in $\H_g$.  
According to the Ansatz of D'Hoker and Phong, we wish to find a modular form, $\Xi^{(g)}[0]$,  of
weight ~$8$ with respect to $\Gamma_g(1,2)$ possessing
the splitting property~$(1)$ and the vanishing trace property~$(2)$:  
\begin{equation}
\label{P1}
\Xi^{(g)}[0]  \pmatrix \tau_1& 0\\ 0& \tau_2\endpmatrix=\Xi^{(g-k)}[0]( \tau_1)\,\Xi^{(k)}[0] (\tau_2)  
\end{equation}
with $\tau_1 \in \Jac_{g-k}$ and $\tau_2 \in \Jac_{k}$
and
\begin{equation}
\label{P2}
\Xi^{(g)}=
\Tr \left(\Xi^{(g)}[0]\right) =
\sum_{\gamma \in   \Gamma_g(1,2)\backslash\Gamma_g } \Xi^{(g)}[0] \vert_8 \gamma =
\sum_{\text{even chars $m$}} \Xi^{(g)}[m]  
\end{equation}
vanishes along $\Jac_g$.\smallskip

Due to the splitting property and the
fact that the accepted solution in genus one,
$ \Xi^{(1)}[0](\tau)=\eta(\tau)^{12}\theta_0(\tau)^4$,
is a cusp form, all of the $\Xi^{(g)}[0]$ should be cusp forms
on the Jacobian locus.  
For $g \le 4$, the $\Xi^{(g)}[0]$ will also be Siegel cusp forms.  
The existence of cusp forms with the properties~$(1)$ and~$(2)$ provides some
vindication of the Ansatz of D'Hoker and Phong. Further vindication would be provided
by the uniqueness of the restriction to the moduli space of curves, the relevant domain for physics. We formulate this as a third condition.  
\begin{equation}
\label{P3}  
\text{Let  }\Xi^{(1)}[0]=\eta^{12}\theta_0^4.  \ \
\text{If } \Xi^{(g)}[0], \,{\tilde \Xi^{(g)}[0]} \in \left[\Gamma_g(1,2),8\right]
\text{ both  satisfy   }
\end{equation}
properties~{$(1)$} and~$(2)$,
then $ {\Xi^{(g)}[0]} = \tilde \Xi^{(g)}[0] $
upon restriction to $\Jac_g$.
\smallskip


The splitting property~$(1)$ means that the chiral superstring measure
is the product measure on the block diagonal.  
The trace property~$(2)$ means that the cosmological constant
vanishes.  The uniqueness property~$(3)$ means that the
desired family of solutions, $\{\Xi^{(g)}[0]\}$, is uniquely determined
on the Jacobian locus by the genus~one solution $\Xi^{(1)}[0]$.  


D'Hoker and Phong gave an expression for $\Xi^{(2)}[0]$ and they proved that is well defined (uniqueness). In \cite{CDPvG} the existence and the uniqueness of  $\Xi^{(3)}[0]$ has been proven.
 The existence of $\Xi^{(4)}[0]$ has been proven in \cite{GR} and \cite{DG}. Moreover in  \cite{CDvG2} a relative uniqueness for $\Xi^{(4)}[0]$ has been  proven. This relative uniqueness is due to two facts.\smallskip
 
 First in genus 4, there is a cusp form of weight 8 with respect to  
 the full integral symplectic group $ \Gamma_4$, the  Schottky form $J$ that vanishes
 along  $\Jac_4$ and hence along the block diagonal period matrices in genus~$4$.  
 If $\Xi^{(4)}[0]$ satisfies properties~$(1)$ and~$(2)$, then so does
 $\Xi^{(4)}[0] + c J$; thus the two desired properties cannot be uniquely
 satisfied by a Siegel modular cusp form in genus four. Accordingly, the
 uniqueness property~$(3)$ only requires the uniqueness of the restriction
 of $\Xi^{(g)}[0]$ to the moduli space of curves and this is consistent with the
 treatment in physics.

 Second, in \cite{CDvG2} the uniqueness of $\Xi^{(4)}[0]$ is proved only among the modular forms
that are polynomials in the second order  theta-constants.   We know, as a consequence of the results in
\cite{Ru} and \cite{Ru2}, that when $g\leq 3$ all modular forms with respect to  $\Gamma_g(1,2)$ are polynomials in the second order  theta-constants. Recently in \cite{OS} two of the authors proved that this statement is false when $g\geq 4$. \smallskip
 
Hence, a priori, in genus 4, there could be other modular forms with the properties
(\ref{P1}) and (\ref{P2})
that cannot be written
as polynomials in the second order theta-constants.  
Our first goal in this article is to prove the uniqueness  property~$(3)$ in genus~$4$
and thus complete the program begun by D'Hoker and Phong through genus~$4$.
The proof is based on the knowledge of the optimal slope of a cusp form in genus~$4$.

Our second goal is to construct $\Xi^{(5)}[0]$ in genus~5.  
In every genus, we classify the linear relations among the eight classes
of self-dual
lattices of rank~{16}. We use this knowledge to construct $\Xi^{(5)}[0]$
as a linear combination of theta series.  We prove the relative uniqueness of
$\Xi^{(5)}[0]$:  it is the only linear combination of theta series in genus~5
satisfying properties~$(1)$ and~$(2)$.  
\smallskip

Using the methods of classical automorphic forms,  
we give independent proofs in all genera $g \le 5$, although  
for motivation we are highly indebted to \cite{CDvG2},
\cite{DHP1}, \cite{DHPa}, \cite{DHPb}, \cite{DHPc}, \cite{GR} and \cite{RSM}.  
To state our main Theorem,
we introduce
$J^{(g)}=\vartheta_{E_8\oplus E_8}^{(g)}-\vartheta_{D_{16}^+}^{(g)}$
for arbitrary genus $g$.
If $g=4$, we get the Schottky form $J$ which we already mentioned.
We prove:  

\begin{thm}
\label{main}
Let $\vartheta^{(g)}$ be the vector of the genus~$g$ theta series of the six
odd self-dual rank~$16$ lattices.  
Let $\{ \Xi_j \}_{j=0}^5 \subseteq \C^6$ be the dual basis to
$\{\Xi^j\}_{j=0}^5 \subseteq \C^6$ for
$\Xi^j=(0,2^j,4^j,8^j,16^j,32^j)$ and $\Xi^0=(1,1,1,1,1,1)$.  
Set $\Theta_j^{(g)}=\Xi_j \cdot \vartheta^{(g)}$.  
For all $g \ge 0$, we have
\begin{equation}
\label{WE}
\vartheta^{(g)} = \Xi^5\, \Theta_5^{(g)}
+ \Xi^4\, \Theta_4^{(g)}+ \Xi^3\, \Theta_3^{(g)}+ \Xi^2\, \Theta_2^{(g)}+ \Xi^1\,  \Theta_1^{(g)}+ \Xi^0\, \Theta_0^{(g)}  .  
\end{equation}
Set $\Xi^{(g)}[0]=\Theta^{(g)}_g-
\frac{17 \cdot 89 \cdot 227}{2^{19} \cdot 3\cdot 5\cdot 7^2\cdot 33}
J^{(g)} \in [\Gamma_g(1,2),8]$.  
For $g \le 4$, the $\Xi^{(g)}[0]$ are cusp forms but $\Xi^{(5)}[0]$ is only a cusp form
on the Jacobian locus.  The family $\{  \Xi^{(g)}[0]  \}$ satisfies both properties~$(1)$
and~$(2)$ of the Ansatz for $g \le 5$ and property~$(3)$ for  $ g \le 4$.  
Also,  $\Xi^{(5)}[0]$ is the unique linear combination of theta series  in
$[\Gamma_5(1,2),8]$ that satisfies both properties~$(1)$
and~$(2)$.  
\end{thm}
And so we leave genus~$5$ where we found genus~$4$,
with the construction of a Siegel modular form that satisfies Ansatz properties~$(1)$
and~$(2)$ but whose uniqueness is only proven among the span of the theta series.
The uniqueness property~$(3)$ for $\Xi^{(5)}[0]$ would follow if we could prove
$\dim [\Gamma_5(1,2),8]=8$.  

\bigskip

We thank the referee for helpful suggestions.

\section{Notation and known results}

For a domain ${\mathbb D }\subseteq \C$,
let $V_g( {\mathbb D } )$ be the $g$-by-$g$ symmetric matrices with coefficients in ${\mathbb D}$.
For ${\mathbb D} \subseteq \R$,
let $\Pg({\mathbb D})^{\text{semi}} \subseteq V_g({\mathbb D})$
be the semidefinite elements and let $\Pg({\mathbb D})$ be the definite elements.  
Let $\H_g$ be the Siegel
upper half~space of genus $g$, i.e. the set of $g$-by-$g$ symmetric complex matrices with positive~definite imaginary part. The symplectic group
$\Gamma_g=\Sp(g,\Z)$ acts on $\H_g$  via
$$
 \pmatrix A&B\\  C&D\endpmatrix\circ\tau:= (A\tau+B)(C\tau+D)^{-1} .
$$
Here we think of elements of $\Gamma_g$ as consisting of four
$g\times g$ blocks.

For $r\in \tfrac12\Z$ and  $\gamma \in \Gamma_g$,  we set
$$(f\vert_r \gamma)(\tau)=\det(C\tau+D)^{-r}f(\gamma\circ\tau)$$
for some choice of square root.  
Let  $\Gamma$ be a subgroup of finite index in $\Gamma_g$, we say that
a holomorphic function $f$ defined on $\H_g$ is  a modular form
of weight $r$ with respect to $\Gamma$  if
$$
\forall\gamma\in\Gamma,\forall\tau\in\H_g, \quad
(f\vert_r \gamma)(\tau)=f(\tau).
$$
and if additionally $f$ is holomorphic at all cusps when $g=1$.
We denote by $[\Gamma, r  ]$ the vector space of such functions.  

For holomorphic $f:\H_g\to \C$
we define
$$\Phi(f)(\tau_1)=\lim_{ \lambda \longrightarrow + \infty}f \pmatrix \tau_1 &0\\ 0&i\lambda\endpmatrix$$
when this limit exists
for all $\tau_1\in\H_{g-1}$. In particular,  this operator maps
$[\Gamma_g, r ] $ to $[\Gamma_{g-1}, r]$ and
$[\Gamma_g(1,2), r ] $ to $[\Gamma_{g-1}(1,2), r]$.  
This operator has a relevant importance in the theory of modular forms,
we refer to \cite{I} or \cite{Fr} for details.  
In the case of the full modular group, a cusp form is a modular form
that is in the kernel of the $\Phi$ operator.
In the case of a subgroup of the modular group, a  modular form is a cusp form if
$$
\forall\gamma\in\Gamma_g, \quad
\Phi( f\vert_r \gamma)=0.
$$
We shall denote by $[\Gamma,\, k]_0$ the subspace of cusp forms.
We shall also use the Witt homomorphism
$
\Psi^{*}_{ij}: \left[\Gamma_{i+j}(1,2),k \right]\to
\left[\Gamma_{i}(1,2),k \right] \otimes  \left[\Gamma_{j}(1,2),k \right]
$
that is the pullback of the map
$\Psi_{ij}:\H_i \times \H_j\to \H_{i+j}$  
defined by $\Psi_{ij}(\tau_1,\tau_2)=\pmatrix \tau_1& 0\\ 0& \tau_2\endpmatrix$.  
For all $f\in [\Gamma_g(1,2),k]$ we have
the following formula for the Fourier coefficients of the image  $\Psi^{*}_{ij}f$:  
\begin{equation}
\label{A2}
a\left(T_1\otimes T_2; \Psi_{ij}^{*}f\right)=
\sum_{T={\tiny\begin{pmatrix} T_1 & * \\ * & T_2 \end{pmatrix}}
\in {\tfrac12 {\cal P}_g(\Z)^{\text{semi}} }} a(T;f).
\end{equation}

\section{The Theta Group}

Before proceeding to the theta group,
we recall the theta functions.
For
$\tau\in\H_g$, $z\in \C^g$ and $\e,\ \de\in \F^g$,
where $\F$ denotes the abelian group $\Z/2\Z=\lbrace 0,1\rbrace$,  
the associated theta function with characteristic $m=[\e, \de]$ is
\begin{align*} \t_m(\tau, z)&=
 \tt\e\de(\tau,z)\\
&=\sum\limits_{n\in\Z^g}e\left( 1/2\cdot
 (n+\e/2)'\tau (n+\e/2)+ (n+\e/2)'( z+\de/2)\right).
\end{align*}
Here we denote by $X'$ the transpose of $X$ and
$e (\ \ )$ stands for $\exp(2\pi i \ \ )$.
As a function of $z$, $\t_m(\tau, z)$ is odd or even depending on whether
the scalar product $\e\cdot\de\in\F$ is equal to 1 or 0,
respectively. Theta-constants are restrictions of theta
functions to $z=0$.
The product of two
theta-constants is a modular form of weight $1 $ with respect to a
subgroup of finite index of $\Gamma_g$.

Now, we discuss the basic properties of the theta group, $\Gamma_g (1,2)$.  
For a $g$-by-$g$ real matrix $X$, we let $X_0\in \R^g$ denote its diagonal.  
We write  $\Gamma_g (n,2n)$  for the subgroup of $\Gamma_g$ defined by
$\gamma\equiv 1_{2g}\,\,{\rm mod} \,n$ and
$$ (AB')_0\equiv  (CD')_0\equiv 0\,\, {\rm mod}\, 2n.$$
Because the theta group is stable under transpose,
we may also use the conditions
$ (A'C)_0\equiv  (B'D)_0\equiv 0\,\, {\rm mod}\, 2n$.

The significance of the theta group $\Gamma_g(1,2)$ is the following:  
$\Gamma_g$ acts on $\F^{2g} $ via
$$
\bmatrix a \\ b \endbmatrix \cdot \smtwomat{A}BCD = \smtwomat{A}BCD' \bmatrix a \\ b \endbmatrix
+ \bmatrix (A'C)_0 \\ (B'D)_0 \endbmatrix.  
$$
If we set $\zeta= \bmatrix a \\ b \endbmatrix$ and
$\gamma=\smtwomat{A}BCD$,
we may write this more compactly
as $\zeta \cdot \gamma = \gamma'\zeta+ \epsilon(\gamma)$,
where $\epsilon(\gamma)= \bmatrix (A'C)_0 \\ (B'D)_0 \endbmatrix$.
The existence of this action follows from the classical transformation of the
theta-constants
$$
\forall \gamma \in \Gamma_g, \quad
\theta[\zeta]{\vert}_{1/2}\gamma \in e\left( \tfrac18\Z\right) \theta[\zeta\cdot\gamma].  
$$
{From} the definitions
it is not difficult to see that
$$
\Gamma_g(1,2)=\Stab_{\Gamma_g}\bmatrix 0 \\ 0 \endbmatrix
\text{ and }  \theta\bmatrix 0 \\ 0 \endbmatrix^8
\in \left[ \Gamma_g(1,2),4 \right].
$$


The group $\Gamma_g$ acts transitively on the even characteristics,
so that the index of $\Gamma_g(1,2)$ in $\Gamma_g$ is
$2^{g-1}(2^g+1)$ and we have the coset decomposition
$$
\Gamma_g= \bigcup_{\text{even $\zeta$}} \Gamma_g(1,2)\, \gamma_{\zeta},
$$
where $\gamma_{\zeta} \in \Gamma_g$ is any element with
$\bmatrix 0 \\ 0 \endbmatrix\cdot \gamma_{\zeta} =\zeta$ in $\F^{2g}$.  
The stabilizers of the nonzero characteristics are given by the conjugate subgroups
$\Stab_{\Gamma_g}[\zeta]=  \gamma_{\zeta}\inv  \Gamma_g(1,2) \gamma_{\zeta}$.  

We let $\Delta_g(\Z)$ denote the subgroup of $\Gamma_g$ with
``$C=0$." The double coset decomposition of $\Gamma_g$ with
respect to $\Gamma_g(1,2)$ and $\Delta_g(\Z)$ is of primary
interest to us. The following Proposition shows that there are just two
double cosets: the $I$-cusp, TIC, and the other cusp, TOC.  
\begin{prop}
\label{B1}
We have the double coset decomposition
$$
\Gamma_g =  \Gamma_g(1,2) \Delta_g(\Z) \cup
\Gamma_g(1,2) \pmatrix I & 0 \\ I & I \endpmatrix \Delta_g(\Z),
$$
where $I=I_g$ denotes the identity matrix of degree $g$.
The first double coset contains $2^g$ single cosets, the second $2^{g-1}(2^g-1)$.  
\end{prop}
{\it Proof.}  
We prove the double coset decomposition  by relating it to
the single coset decomposition, so we first give representatives
$\gamma_{\zeta}$ for the single coset decomposition.  
For even $\zeta=\bmatrix a \\ b \endbmatrix$, we may choose $\gamma_{\zeta}$ as
$$
\gamma_{\zeta}=
\pmatrix I & 0 \\ \diag(a) & I \endpmatrix
\pmatrix I & S \\ 0 & I \endpmatrix=
\pmatrix I & S \\ \diag(a) & I+\diag(a)S \endpmatrix
$$
with any $S\in V_g(\Z)$ such that $Sa+S_0=b$.  
When $a=0$, we just take $S=\diag(b)$.  
When $a \ne 0$, one way to get $S$ is to
take $S= \beta b'+b \beta'$ for any {\it odd\/} characteristic
$\bmatrix a \\ \beta \endbmatrix$.    

We will show that
the right action of $\Delta_g(\Z)$ on the set of the even characteristics has two orbits:
$\bmatrix 0 \\ * \endbmatrix$ with $2^g$ elements and
$\bmatrix \ne 0 \\ * \endbmatrix$ with $2^{g-1}(2^g-1)$ elements.  
For $\zeta=\bmatrix 0 \\ b \endbmatrix$, we have
$ \gamma_{\zeta }\in   \Gamma_g(1,2)  \Delta_g(\Z) $ because
$
\gamma_{\zeta }= \smtwomat{I}{\diag(b)}{0}{I}
$.  
For $\zeta=\bmatrix a \\ b \endbmatrix$, with $a \ne 0$, we have
$ \gamma_{\zeta }\in   \Gamma_g(1,2) \smtwomat{I}{0}{I}{I} \Delta_g(\Z) $ because
for any $U \in \GL_g(\Z)$ we have
\begin{align*}
\gamma_{\zeta }  &=\pmatrix {I} &{0} \\ {\diag(a)} &{I} \endpmatrix
\pmatrix {I} &{\beta b'+b \beta'} \\ {0} &{I} \endpmatrix  \\
&=
\pmatrix {U\inv} &{0} \\ {\diag(a)U\inv-U'} &{U'} \endpmatrix
\pmatrix {I} &{0} \\ {I} &{I} \endpmatrix
\pmatrix {U} &{0} \\ {0} & (U')^{-1} \endpmatrix
\pmatrix {I} &{\beta b'+b \beta'} \\ {0} &{I} \endpmatrix.    
\end{align*}
To complete the proof we show that
$\smtwomat{U\inv}{0}{\diag(a)U\inv-U'}{U'}\in \Gamma_g(1,2)$
for  some choice of $U $.  Choose $U$ so that $U'I_0 \equiv a \mod 2$.  
The defining conditions,  
$ (AB')_0\equiv  (CD')_0\equiv 0\,\, {\rm mod}\, 2$,
are then satisfied if we note that $(U'U)_0 \equiv U'I_0 \mod 2$.  \qed

The choice of
$\gamma_{\zeta}=\smtwomat{I}0II$ and of $\zeta=\bmatrix I_0 \\ 0 \endbmatrix$
as a representative for the nontrivial double coset
is somewhat arbitrary. We have chosen the even characteristic
$\bmatrix I_0 \\ 0 \endbmatrix$
because its upper and lower parts are invariant
under permutations and the representative $\gamma_{\zeta}$
because it is a direct sum of $\SL_2(\Z)$ matrices.  

We now study the Fourier expansions at these two cusps.  
For $S\in V_g(\R)$, let $t(S)=\smtwomat{I}{S}{0}{I} \in \Sp_g(\R)$.  
For $U\in \GL_g(\R)$, let
$u(U)=\smtwomat{U}{0}{0}{(U')^{-1}} \in \Sp_g(\R)$.
For $A,B \in V_g(\C)$, set $\<A,B \>= \tr(AB)$.  

\medskip

{\bf{The $I$-Cusp.} }  
The translation lattice for the $I$-cusp (TIC) is
$$
\{ S \in V_g(\Z): t(S) \in \Gamma_g(1,2) \}=
\{ S \in V_g(\Z): \text{ $S$ even} \}.
$$
The dual lattice
with respect to $\< \ \ ,\ \  \>$ is $\frac12 V_g(\Z)$.  
If we multiply the above translation lattice by $1/2$, we get the lattice $\Xg$
consisting of ``half-integral" matrices, which is nothing else but the dual lattice of $V_g(\Z)$.
The similarity lattice for TIC is
$$
\{ U \in \GL_g(\Z): u(U) \in \Gamma_g(1,2) \}=
\GL_g(\Z).
$$
Therefore, an $f \in \left[ \Gamma_g(1,2),k \right]_0$ has a
Fourier expansion at TIC
$$
f(\tau)= \sum_T a(T) e\left( \<\tau, T \>\right)
$$
where $T \in \frac12 {\cal P}_g(\Z)$ runs over integral forms multiplied
by $\frac12$ and, for all $ U \in \GL_g(\Z)$,
$a(U' T U)= \det(U)^k \,a(T)$.  

\medskip

{\bf{The Other Cusp.}  }
Let $M = \smtwomat{I}0II$ and let  $\Gamma_g(1,2)^M$ denote
$M\inv \Gamma_g(1,2) M$.  
The translation lattice ${\cal L\/} $ for the other cusp (TOC) is
given by  $S$ such that $t(S)$ stabilizes $\bmatrix I_0 \\ 0 \endbmatrix$.  

\begin{align*}
{\cal L\/}  &=
\{ S \in V_g(\Z): t(S) \in \Gamma_g(1,2)^M \} \\
&=
\{ S \in V_g(\Z): \text{ $SI_0+S_0 \equiv 0\mod 2$} \}.  
\end{align*}

\begin{lm}
\label{B2}
Let $E_{ij}\in  V_g(\Z)$ be the matrix with $1$ in the $(i,j)$ and $(j,i)$ spots
and zeroes elsewhere.  
The  lattice ${\cal L\/} $ contains
$2 V_g(\Z)$, $\diag(\Z^g)$ and  
$E_{ij}+E_{jk}+E_{ki}$ for all distinct triples
$i,j,k \in \{1,2,\dots, g\}$.  
\end{lm}
{\it Proof.}
This is an easy computation from the defining condition:
$SI_0+S_0 \equiv 0\mod 2$.  \qed

Actually, the above elements span  ${\cal L\/} $,
although the last group $E_{ij}+E_{jk}+E_{ki}$ is linearly dependent over $\F$ for $g > 3$.  
Furthermore, the indices are
$[V_g(\Z):{\cal L\/}]=2^{g-1}$  and
$[{\cal L\/}:2 V_g(\Z)]=2^{1+\binom{g}{2}}$ but we need none of these facts.

\begin{df}
\label{B3}
Elements of the lattice $4 {\Cal L\/}^{*} $ are called {\bf very even\/}.      
\end{df}
{From} $ 2 V_g(\Z) \subseteq {\Cal L\/} \subseteq V_g(\Z)$ we see that the dual
lattice, ${\Cal L\/}^{*}$, satisfies
$\tfrac12 \Xg \supseteq {\Cal L\/}^{*} \supseteq  \Xg $, so that the elements of
$4 {\Cal L\/}^{*} $ are in fact even.  
Since $\diag(\Z^g) \subseteq {\Cal L\/}$, the diagonals of
${\Cal L\/}^{*}$ are integral and the diagonals of
very even forms are multiples of~$4$.

The similarity group ${\Cal G\/}$  for TOC is given by the $U$ such that
$u(U)$ stabilizes $\bmatrix I_0 \\ 0 \endbmatrix$:
\begin{align*}
{\Cal G\/}  &=
\{ U \in \GL_g(\Z): u(U) \in \Gamma_g(1,2)^M \} \\
&=
\{ U \in \GL_g(\Z): \text{ $U'I_0 \equiv I_0\mod 2$ } \}.  
\end{align*}
A useful comment here is that ${\Cal G\/}$ contains all permutation
matrices and all diagonal sign changes.  

Therefore, an $f \in \left[ \Gamma_g(1,2),k \right]_0$ has a
Fourier expansion at TOC
$$
(f|M)(\tau)= \sum_T b(T) e\left( \<\tau, T \>\right)
$$
where $T \in {\Cal L\/}^{*}\cap \Pg(\Q)$ runs over very even forms multiplied
by $\frac14$ and, for all $ U \in {\Cal G\/}$,
$b(U' T U)= \det(U)^k\, b(T)$.  

\begin{prop}
\label{B4}
Let $f \in \left[ \Gamma_g(1,2), k \right] $.  
The form $f$ is a cusp form if and only if
$\Phi(f)=0$ and  $\Phi(f|M)=0$.  

Furthermore, let the coset representatives for $\Gamma_g(1,2)
\backslash \Gamma_g$ be written
$$
\gamma_{\zeta}=
\pmatrix I & 0 \\ \diag(a) & I \endpmatrix
\pmatrix I & S \\ 0 & I \endpmatrix=
M_a \, t(S).  
$$
For $S\in V_g(\Z)$, let $\pi(S) \in V_{g-1}(\Z)$
be the lower $(g-1)$-by-$(g-1)$ block.  \newline
For $a =(0; a')$, we have
$ \Phi( f | M_a \, t(S))= \Phi( f )  |  M_{a'} \, t(\pi(S))$.  \newline
For $a =I_0 + (0; c')$, we have
$ \Phi( f | M_a \, t(S))=\Phi(f | M) | M_{c'}\, t(\pi(S))$.  
\end{prop}
{\it Proof.}  
The form $f$ is a cusp form if and only if  $\Phi(f|\gamma_{\zeta})=0$ for
every $\gamma_{\zeta}$ in the single coset decomposition of
$\Gamma_g(1,2) \backslash \Gamma_g$.  
We rely on two properties of the $\Phi$ map.  
First, for any  $F \in \left[ \Gamma, k \right] $, we have
$\Phi(F | t(S))=\Phi(F) | t( \pi(S) )$.  
Second, let $I_2 \in \Gamma_1$ be the identity matrix;
for $\gamma' \in \Gamma_{g-1}$, we have
$\Phi(F | I_2 \oplus \gamma' )=\Phi(F) | \gamma' $.  
Here we understand that
\[
\begin{pmatrix} a & b \\ c & d \end{pmatrix}
\oplus \begin{pmatrix} A & B \\ C & D \end{pmatrix}
=\begin{pmatrix} a & 0 & b & 0 \\
0 & A & 0 & B \\
c & 0 & d & 0 \\
0 & C & 0 & D \end{pmatrix}.
\]
 
 There are two cases depending upon $a$.  
 If $a$ begins with $0$, so that $a =(0; a')$,
 then $M_a= I_2 \oplus M_{a'}$ and
\begin{multline*}
\Phi( f | M_a \, t(S))= \\
\Phi( f | M_a)  |  t(\pi(S))=
\Phi( f | I_2 \oplus M_{a'} )  |  t(\pi(S))=
 \Phi( f )  |  M_{a'} \, t(\pi(S)).
 \end{multline*}
 If $a$ begins with $0$, define $c= a - I_0$;
 then    $c$ begins with $0$, so that $c =(0; c')$.  
 Note that $M=M_{I_0}$ and that $M_a=M M_c$.  
 We have
   $$
 \Phi( f | M_a \, t(S))=
  \Phi( f | M M_c)  |  t(\pi(S))=
\Phi(f | M) | M_{c'}\, t(\pi(S)).  
$$
Thus,  the cases $\Phi(f|\gamma_{\zeta})=0$
follow from $\Phi(f)=0$ and  $\Phi(f|M)=0$.  \qed

\section{{\bf Lemmata}}

Let $\Gamma \subseteq \Gamma_g$ be a subgroup
of finite index~$I=[ \Gamma_g:\Gamma]$. The Norm map is defined
as
\begin{align*}
\Norm: \left[ \Gamma, k \right] &\to \left[ \Gamma_g, Ik \right]  \\
f &\mapsto \prod_{\gamma\in \Gamma \backslash \Gamma_g} f \vert_k \gamma  .
\end{align*}
This map naturally induces the map between the associated projective spaces
and we use the same notation $\Norm$ again.
The next Lemma shows that if we have a nontrivial subspace~$S$ of
Siegel cusp forms, all of whose elements have a
norm that is a multiple of a fixed form,
then the dimension of $S$ is one.
\begin{lm}
\label{Q0}
Let $\Gamma \subseteq \Gamma_g$ be a subgroup
of finite index~$I$.  
Let $S \subseteq \left[ \Gamma, k \right]_0$ be a subspace.  
If $ \Norm: \Pj\left( S \right) \to \Pj\left( \left[ \Gamma_g, Ik \right]_0 \right)$
has an image consisting of precisely one point, then $\dim S =1$.  
\end{lm}
{\it{Proof.}}
Take $A,B\in S \setminus \{0\}$.  
We will show that $B/A$ is constant as a meromorphic function on $\H_g$,
and thus conclude that $\dim S=1$.  
Let $\xi=\Norm(A) \in  \left[ \Gamma_g, Ik \right]_0$ and note that $\xi \ne 0$
since $A \ne 0$.  By assumption, for all $x \in \C$ there exists a $c \in \C$
such that $\Norm(Ax+B)=c\,\xi$.  We will show that $c$ is a polynomial in $x$.  
We can evaluate $c$ by picking $\tau_0 \in \H_g$ with
$\xi(\tau_0) \ne 0$; then
\begin{align*}
c=  &\frac{\Norm(Ax+B)(\tau_0)}{\xi(\tau_0)}
= \prod_{\gamma}^I \left(     x+ \frac{(B\vert \gamma)(\tau_0)}{(A\vert \gamma)(\tau_0)} \right).   \\
\end{align*}
However, since $c = \Norm(Ax+B)/\xi$ as well, we have
\begin{align*}
&\prod_{\gamma}^I \left(     x+ \frac{(B\vert \gamma)}{(A\vert \gamma)} \right)=
\prod_{\gamma}^I \left(     x+ \frac{(B\vert \gamma)(\tau_0)}{(A\vert \gamma)(\tau_0)} \right).  
\end{align*}
Letting $x = -B/A$, we see that $B/A$ is a meromorphic function
with a discrete image and hence is a constant.  
\qed

\begin{df}
\label{E1}
A function $\phi:\PsR\to\R_{\ge 0}$
is called
{\it{type one}} if
\item{1.}   For all $s \in \PR$, $\phi(s) > 0$,  
\item{2.}   for all
$\lambda\in\R_{\ge0}$ and $s\in\PsR$, $\phi(\lambda s)=\lambda\phi(s)$,
\item{3.}   for all $s_1,s_2\in\PsR$,  $\phi(s_1+s_2)\ge\phi(s_1)+\phi(s_2)$.
\end{df}
Type one functions are continuous on $\PR$ and
respect the partial order on $\PsR$.    
We will need the acquaintance of two type one functions:
the Minimum function $m(s) = \inf_{u\in\Z^g\backslash\{0\}}u'su$
and its convex dual $w(s) = \inf_{u\in\PR}{\<u,s\>}/{m(u)}$, the dyadic trace.  

Recall the definition of the {\it slope\/.}
\begin{align*}
\text{For $ f \in \left[ \Gamma, k \right]_0$,  let }
&\supp(f)=\{ T \in \Pg(\Q): a(T;f) \ne 0 \};  \\  
m(f) &= \min m\left( \supp(f) \right) \text{ and }
\slope(f)=\frac{k}{m(f)}.  
\end{align*}
We know that the minimal slope on
$\A_g=\Gamma_g\backslash \H_g$ is: $12$, $10$, $9$, $8$ for
$g=1,2,3,4$, respectively, cf. \cite{HM} and in particular the
Corollary to Theorem 3.
It is  rather easy to check that
any Siegel modular form that attains this slope is a power of
$\Delta\in [\Gamma_1,12]_0$,
$X_{10}\in [\Gamma_2,10]_0$,
$X_{18}\in [\Gamma_3,18]_0$ or
$J\in [\Gamma_4,8]_0$ in $g=1,2,3$ or $4$, respectively.

In the next Lemma,  we extend this result to subgroups $\Gamma$ of index $I$ by
looking at the slope of the average vanishing, $k/ \mu(f)$, where
$$
\mu(f)= \frac1{I} \sum_{\gamma \in \Gamma\backslash\Gamma_g}^I m\left( f \vert \gamma\right).  
$$

\begin{lm}
\label{Q1}
Let $\Gamma \subseteq \Gamma_g$ be a subgroup
of finite index~$I$.
Let $f \in \left[ \Gamma, k \right]_0$.  
If
$k/ \mu(f)$  is the optimal value: $12$, $10$, $9$, $8$ for
$g=1,2,3,4$, respectively,
then $\Norm(f) \in  \left[ \Gamma_g, Ik \right]_0$ is a constant multiple
of a power of
$\Delta,\ X_{10},\ X_{18}, \ J$, respectively.
\end{lm}
{\it{Proof.}}
The slope of the level one $\Norm(f)$ is $Ik/m \left( \Norm(f) \right)$.  
Therefore, it suffices to prove that
$ m \left( \Norm(f) \right) \ge I \mu(f)$.
Note that $m$ is a type one function satisfying $m(s_1+s_2) \ge m(s_1) + m(s_2)$.  
 We have
\begin{align*}
m \left( \Norm(f) \right)
&=  m \left(  \prod_{\gamma\in \Gamma \backslash \Gamma_g} f \vert \gamma   \right)
= \min  m \left( \supp \left(  \prod_{\gamma}^I f \vert \gamma   \right)  \right) \\
&\ge \min  m \left(   \sum_{\gamma}^I \supp(f \vert \gamma)   \right)  
\ge  \min  \sum_{\gamma}^I m\left(\supp(f \vert \gamma) \right) \\
&= \sum_{\gamma}^I  \min m(\supp(f \vert \gamma) )
=  \sum_{\gamma}^I m(f \vert \gamma)
= I \mu(f).    \qed
\end{align*}

The dyadic trace, defined as
$w(s) = \inf_{u\in\PR}{\<u,s\>}/{m(u)}$,  
also has a characterization as a supremum \cite{PY00}.  
A {\it dyadic representation\/} of a form $T \in {\cal P}_g(\Q)$ is
given by $\a_i>0$ and
$v_i \in \Z^g\backslash \{0\}$ such that
$T=\sum_{i} \a_iv_iv_i'$.  Since
$\<T,u\>= \sum_{i} \a_i \<v_iv_i',u\> \ge \sum_{i} \a_i m(u)$,
it follows from the definition of the dyadic trace that
$w(T) \ge \sum_{i} \a_i$ for any dyadic representation.  
In fact, we have
$$
w(T) = \sup \{ \sum_{i} \a_i: \text{ over all dyadic representations $\sum_{i} \a_iv_iv_i' $ of $T$}  \}
$$
and that this supremum is attained by a particular dyadic representation.  
The dyadic trace is a useful tool in the geometry of numbers.  

\begin{lm}
\label{Q2}
Let $T \in \Pg(\Q)$ with $g=g_1+g_2$,
$T=\smtwomat{T_1}{W}{W'}{T_2}$
for $T_1 \in \Pgone(\Q)$,
$T_2 \in \Pgtwo(\Q)$,
$W \in \operatorname{Mat}_{g_1 \times g_2}(\Q)$.  
We have
$$
w(T) \le w(T_1) + w(T_2).  
$$
Furthermore, we have equality if and only if $W=0$.  
\end{lm}
{\it{Proof.}}  
Let $T=\sum_{i} \a_iv_iv_i'$ with $\a_i >0$ and
$v_i \in \Z^g\setminus\{0\}$ be a dyadic
representation of $T$ that attains the dyadic trace:
$w(T) =\sum_i \a_i$.  
We use $\pi_1(v),\ \pi_2(v)$ to denote
the first $g_1$ and the last $g_2$ coordinates of $v$.
For $j=1,2$
 $$
 T_j=\sum_{i: \pi_j(v_i) \ne 0} \a_i \pi_j(v_i)\pi_j(v_i)'
 $$
is a dyadic representation of $T_j$ so that
$w(T_j) \ge \sum_{i: \pi_j(v_i) \ne 0} \a_i $.  
Therefore we have
\begin{align*}
w(T_1) + w(T_2)  & \ge
\sum_i \a_i +  \sum_{i: \pi_1(v_i) \ne 0 \text{ and } \pi_2(v_i) \ne 0} \a_i
\ge \sum_i \a_i =w(T).  \\
\end{align*}
This is the first assertion.  
Equality can be attained only if the second sum above is empty.  
In this case we have
\begin{align*}
T &=\sum_{i} \a_iv_iv_i' =
\sum_{i} \a_i  \begin{pmatrix}\pi_1 v_i   \\  \pi_2 v_i \end{pmatrix}
\begin{pmatrix}\pi_1 v_i'  & \pi_2 v_i' \end{pmatrix} \\
& =  \sum_{i: \pi_1(v_i)=0} \a_i  \begin{pmatrix}\pi_1 v_i   \\  \pi_2 v_i \end{pmatrix}
\begin{pmatrix}\pi_1 v_i'  & \pi_2 v_i' \end{pmatrix}
+ \sum_{i: \pi_2(v_i)=0} \a_i  \begin{pmatrix}\pi_1 v_i   \\  \pi_2 v_i \end{pmatrix}
 \begin{pmatrix}\pi_1 v_i'  & \pi_2 v_i' \end{pmatrix} \\
&= \smtwomat000{T_2} + \smtwomat{T_1}000=\smtwomat{T_1}00{T_2} .  
\end{align*}
This shows $W=0$ and  so
$ w(T) = w(T_1) + w(T_2)$ indeed holds.  \qed

\begin{lm}
\label{Q3}
Let  $ f \in  \left[ \Gamma_g(1,2), k \right]$ for $g \ge 2$.  
If $ f \in   \ker \Psi_{1,g-1}^{*}$, then  we have $m(f) \ge 1$.  
\end{lm}
{\it Proof.}  
Recall that the Fourier expansion of $f$ is
$$
f(\tau)= \sum_{T} a(\tfrac12 T;f) e \left( \<\tau, \tfrac12 T \> \right)
$$
with $T \in \Pg(\Z)$ since $ \ker \Psi_{1,g-1}^{*} \subseteq  \left[ \Gamma_g(1,2), k \right]_0$.  
For $ \tfrac12 T \in \supp(f)$ we will show that $m(T) \ge 2$.  
It suffices to prove that $a\left(  \tfrac12(1 \oplus T_0);f\right)=0$ for all
$T_0\in {\cal P}_{g-1}(\Z)$ because, if $m(T)=1$,
then $T$ is $\GL_g(\Z)$-equivalent to $1 \oplus T_0$.  

We will show that $a\left(  \tfrac12(1 \oplus T_0);f\right)=0$ by
induction on $w(T_0)$.  The base case of the induction is satisfied
because $f$ is a cusp form.  
Since $\Psi_{1,g-1}^{*}f=0$, its $ \tfrac12 1 \otimes \tfrac12 T_0 $ Fourier coefficient is also $0$
and
\begin{equation}
\label{F7}
0{= }\sum_{v \in \Z^{g-1}}
a\left(  \tfrac12\smtwomat{1}{v}{v'}{T_0}\right)
{=} a\left(  \tfrac12\smtwomat{1}{0}{0}{T_0}\right){+}
\sum_{v \ne 0}
a\left(  \tfrac12\smtwomat{1}{v}{v'}{T_0}\right)
\end{equation}
All the indices $\smtwomat{1}{v}{v'}{T_0}$ are
$\GL_g(\Z)$-equivalent to $1 \oplus T_v$ for some
$T_v\in {\cal P}_{g-1}(\Z)$.  By Lemma~\ref{Q2} we have
\begin{align*}
1+w(T_v) =
w \smtwomat{1}{0}{0}{T_v}   =    
w \smtwomat{1}{v}{v'}{T_0}   &< 1+w(T_0) \text{ for $v \ne 0$, }   \\
\end{align*}
so that $w(T_v) < w(T_0)$ for   $v \ne 0$.  
By the induction hypothesis, we have
$$
a\left(  \tfrac12\smtwomat{1}{v}{v'}{T_0}\right)=
a\left(  \tfrac12\smtwomat{1}{0}{0}{T_v}\right)
a\left(  \tfrac12(1 \oplus T_{v});f\right)
=0
$$
for   $v \ne 0$ and so
$a\left(  \tfrac12\smtwomat{1}{0}{0}{T_0}\right)=0$
by equation~(\ref{F7}) as well.  \qed

\section{Dimensions}

This next Theorem is a consequence of the work of Igusa  and
the facts previously discussed in genus~three,
cf. \cite{Ru2}, \cite{CDPvG}.  

\begin{thm}
\label{T1}
For $g=1,2,3$, we have $\dim \left[ \Gamma_g(1,2), 8 \right]_0 =1$.  
\end{thm}
{\it Proof.}  
Recall that  these spaces are nonempty, containing nonzero elements $\vartheta^{(g)}_{1_8,\, P_4}$\
(See the appendix for this function).
The Minimum function  $m$ is a $\GL_g(\Z)$-class function, so that  
$m( f | \gamma)$ only depends upon the double coset
$\Gamma_g(1,2) \gamma \Delta_g(\Z)$.  Hence,
for nontrivial $f \in  \left[ \Gamma_g(1,2), 8 \right]_0$ we have, by
the double coset decomposition of
Proposition~\ref{B1},  (set $I= 2^{g-1}(2^g+1)$)  
$$
\mu(f)=
\frac1{I} \sum_{ \gamma \in \Gamma_g(1,2)\backslash\Gamma_g}^I m\left( f \vert \gamma\right)=
\dfrac{2}{ 2^g+1 }\,m(f) +
\dfrac{ 2^g-1 }{ 2^g+1 }\,m(f|\,\smtwomat{I}0II) .  
\label{C6}
$$
The indices at TIC consist of $\tfrac12$ times integral forms and so
$m(f ) \ge \tfrac12$.  
The indices at TOC consist of $\tfrac14$ times very even forms and so
$m(f|\,\smtwomat{I}0II) \ge 1$.  
Therefore we have
$$
\mu(f) \ge
\dfrac{2}{ 2^g+1 }\,\left( \frac12 \right) +
\dfrac{ 2^g-1 }{ 2^g+1 }\,\left( 1 \right) = \dfrac{2^g}{2^g+1}
\label{C7}
$$
and we can make the following table of the maximum  slope
of the average vanishing
for nontrivial elements of $ \left[ \Gamma_g(1,2), 8 \right]_0$:  

\medskip

\centerline{{\bf  Table 1. \/}  Maximum slope for average vanishing.  }
\begin{align*}
g \hskip2.4in  & k/\mu(f)  \\
1\hskip2.1in \frac{8}{2/3}  &= 12  \\
2\hskip2.1in \frac{8}{4/5}  &= 10  \\
3\hskip2.1in \frac{8}{8/9}  &= 9.  
\end{align*}
By Lemma~\ref{Q1}, we know that $\Norm(f)$ is
some multiple of
$\Delta^2\in [\Gamma_1,24]$,
$X_{10}^8\in [\Gamma_2,80]$,
$X_{18}^{16}\in [\Gamma_3,288]$,  
respectively.  

The image of
$ \Norm:  \Pj\left( \left[ \Gamma_g(1,2), 8 \right]_0 \right)
\to \Pj\left( \left[ \Gamma_g, * \right]_0 \right)$
consists of one point in these cases so that
$\dim \left[ \Gamma_g(1,2), 8 \right]_0 =1$ by Lemma~\ref{Q0}.  \qed

\begin{thm}
\label{T2}
For $g=2,3$, we have
$
\left[ \Gamma_g(1,2), 8 \right] \cap \ker \Psi_{1,g-1}^{*}=\{0\}  
$.  
For $g=4$, we have
$
\left[ \Gamma_4(1,2), 8 \right] \cap \ker \Psi_{1,3}^{*}=\C J
$.  
\end{thm}
{\it Proof.}  
By Lemma~\ref{Q3},
an $f \in  \left[ \Gamma_g(1,2), 8 \right] \cap \ker \Psi_{1,g-1}^{*}$
has $m(f) \ge 1$.  
Therefore, as in the proof of the previous Theorem,
$$
\mu(f) \ge
\dfrac{2}{ 2^g+1 }\,\left( 1 \right) +
\dfrac{ 2^g-1 }{ 2^g+1 }\,\left( 1 \right) =1
$$
and $f$ has $8/\mu(f)$  at most $8$.  
Hence $f=0$ in $g=2,3$ and $\Norm(f)$
is a multiple of $J^{136}$ in $g=4$ by Lemma~\ref{Q1}.  

In $g=4$ therefore, the image of the map
$$ \Norm:  \Pj\left( \left[ \Gamma_4(1,2), 8 \right]  \cap \ker \Psi_{1,3}^{*}   \right)
\to \Pj\left( \left[ \Gamma_4, 8\cdot 136 \right]_0 \right)$$
has just one point.  
Therefore,  we have
$\dim  \left[ \Gamma_4(1,2), 8 \right]   \cap \ker\Psi_{1,3}^{*}  =1$
by Lemma~\ref{Q0}, and necessarily
$ \left[ \Gamma_4(1,2), 8 \right]  \cap \ker \Psi_{1,3}^{*} =\C J$.  \qed

 \section{Linear Relations among Theta series}
A more general way to construct modular forms is to use theta series,
in particular if $L$ is a self-dual lattice of rank $m$,
with $8$ dividing $m$,
then
we have the associated quadratic form $S$ and the theta series
\begin{equation}
\label{T6}
\vartheta_{L}^{(g)}(\tau) =\sum_{X\in \Z^{m,\,g}}  
e(1/2\cdot {\rm tr}(S[X]\tau))\notag
\end{equation}
is a modular form of weight $m/2$ relative to $ \Gamma_g(1,2)$.
We let $[ \Gamma_g(1,2),k]^{\vartheta}$ denote the subspace spanned by
theta series
of self-dual lattices of rank~{$2k$}.  

There are eight self-dual lattices of rank~{16}, 
two even lattices and six odd.  
The theta series are elements of $[ \Gamma_g,8]$ and
$[ \Gamma_g(1,2),8]$, respectively.  
In this section we find all the linear relations among these theta series for
every genus.  We give two applications.  
First, we derive the results $\dim [ \Gamma_4(1,2),8]=7$ and
$\dim [ \Gamma_4(1,2),8]_0=2$.  
Second, we push the physicists' program to success in genus five and prove our
main Theorem~\ref{main}.

The eight self-dual lattices of rank~{16} are all found by Kneser's gluing method.  
We use the notation in  \cite{CS}, Table 16.7. For
$4|n$, $D_n^{+}= D_n \cup ([1]+D_n)$ is unimodular;  
$D_8^{+}$ is commonly denoted by $E_8$.  The two even lattices are given by
$E_8 \oplus E_8$ and $D_{16}^{+}$. An odd lattice is given by
$\Z^4 \oplus D_{12}^{+}$ and another by
\begin{align*}
(D_8 \oplus D_8)^{+} &=
D_8 \oplus D_8 \cup
([1]\times [2]+D_8 \oplus D_8)  \\
\cup\, &([2]\times [1]+D_8 \oplus D_8) \cup
([3]\times [3]+D_8 \oplus D_8) ,  \text{ where } \\
[1]\times [2]  &=
[\frac12, \frac12, \frac12, \frac12, \frac12, \frac12, \frac12, \frac12;
0, 0, 0, 0, 0, 0, 0, 1],  \\
[2]\times [1]  &=
[0, 0, 0, 0, 0, 0, 0, 1;\frac12, \frac12, \frac12, \frac12, \frac12, \frac12, \frac12, \frac12], \\
[3]\times [3]  &=
[\frac12, \frac12, \frac12, \frac12, \frac12, \frac12, \frac12, -\frac12 ;
\frac12, \frac12, \frac12, \frac12, \frac12, \frac12, \frac12, -\frac12].    
\end{align*}
Also, we have $\Z \oplus A_{15}^{+}$ for
$$
A_{15}^{+} =
A_{15}  \cup
([4] +A_{15} ) \cup  
([8] +A_{15} ) \cup
([12] +A_{15} ) ,
$$
where $[i]$, for $i+j=n+1$, means
$j$ coordinates of $\frac{i}{n+1}$ followed by $i$ coordinates of $\frac{-j}{n+1}$.  
Finally, we have $\Z ^2\oplus (E_7 \oplus E_7)^{+}$ for
\begin{align*}
(E_7 \oplus E_7)^{+}   &=
(E_7 \oplus E_7)   \cup
([1] \times [1]+(E_7 \oplus E_7)) , \text{ where } \\
[1]\times [1]  &=
[\frac14, \frac14, \frac14, \frac14, \frac14, \frac14, -\frac34, -\frac34;
\frac14, \frac14, \frac14, \frac14, \frac14, \frac14, -\frac34, -\frac34].  
\end{align*}

The following Table gives basic information about the theta series
$\vartheta^{(g)}_i=\vartheta^{(g)}_{\Lambda_i}$
of
these lattices, labeled as $\Lambda_i$ for $i=0,1,\dots,7$.  
Here, $\tau_i$ is the number of vectors of norm one.  
It is the coefficient of $q^{1/2}$ in the genus one Fourier expansion $\vartheta_i^{(1)}$,
whose leading term is also given.  The space
$\left[ \Gamma_1(1,2),8\right]$ is spanned by
$\Xi^{(1)}[0]$, $\vartheta^{(1)}_0$ and $\vartheta^{(1)}_6$ and the coefficients
$\tau_i$, $b_i$, $c_i$, of the linear relation
$\vartheta_i^{(1)}=\tau_i\,\Xi^{(1)}[0]+b_i\,\vartheta_0^{(1)}+c_i\,\vartheta_6^{(1)}$
are also given.  \smallskip

\centerline{{\bf Table 2.}  The eight self-dual lattices of rank 16. }~
\begin{alignat*}6
&i  &  &\Lambda_i  & &\tau_i &   &b_i \quad& c_i \qquad&  &\vartheta_i^{(1)}-1\hskip1.0in\\
&0\quad  &  &(D_8 \oplus D_8)^{+}\quad  &  &0\quad &  &1& 0 \qquad&  &224q^1+4096q^{3/2} \\
&1\quad  &  &\Z\oplus A_{15}^{+}\quad  &  &2\quad &   &1 & 0 \qquad&  &2q^{1/2}+240q^1+4120q^{3/2} \\
&2\quad  &  &\Z^2\oplus (E_7\oplus E_7)^{+} \,\,\, &  &4\quad &  &1 & 0\qquad&  &4q^{1/2}+256q^1+4144q^{3/2} \\
&3\quad  &  &\Z^4\oplus D_{12}^{+}\quad  &  &8\quad &   &1 & 0 \qquad&  &8q^{1/2}+288q^1+4192q^{3/2} \\
&4\quad  &  &\Z^8\oplus E_8\quad  &  &16\quad &   &1 & 0 \qquad&  &16q^{1/2}+352q^1+4288q^{3/2} \\
&5\quad  &  &\Z^{16} \quad  &  &32\quad &   &1 & 0 \qquad&  &32q^{1/2}+480q^1+4480q^{3/2} \\
&6\quad  &  &E_8 \oplus E_8 \quad  &  &0\quad &   &0 & 1 \qquad&  &480q^1 \hskip0.8in \\
&7\quad  &  &D_{16}^{+}   \quad  &  &0\quad &   &0 & 1 \qquad&  &480q^1 \hskip0.8in
\end{alignat*}

We recall the general set up of \cite{er} for organizing linear relations among
theta series.  
Let $\{L_i\}_{i=1}^h$ be a set of self-dual lattices of
rank~$2k$.  Let $V=\C^h$ and define
$\Theta^{(g)}: V \to  \left[ \Gamma_g(1,2), k \right]$ by
$\Theta^{(g)}(r)=\sum_{i=1}^h r_i \vartheta_{L_i}^{(g)}$
and the convention $ \vartheta_{L_i}^{(0)}=1$.  
We have a decreasing filtration $V_g=\ker \Theta^{(g)}$ of $V$ and the dual
increasing filtration $W_g =(V_g)^{\perp}$ of the dual space $V^*$.  
Composing the canonical isomorphism
$ \Theta^{(g)}( V )\cong V/V_g$ with the noncanonical
$V/V_g \cong W_g$, we have
$\dim [\Gamma_g(1,2),k]^\vartheta = \dim W_g$.

The Fourier coefficients of the theta series provide natural
elements of $W_g$.  
For $T \in \tfrac12 \Pg(\Z)$ and
$\vartheta_{L_i}^{(g)}(\tau)= \sum_T a_i^{(g)}(T) e\left( \<T,\tau\>\right)$,
we define $w_g(T) \in V^*$ by $w_g(T)_i=a_i^{(g)}(T)$.  
The $w_g(T)$ for $T \in \tfrac12 \Pg(\Z)$ span $W_g$
and linear relations may be presented by giving a basis of $W_g$
in terms of the $w_g(T)$
or linear combinations thereof.  
We also define the component-wise multiplication on $V^*$ because
this multiplication respects the $W_g$-filtration: $W_iW_j \subseteq W_{i+j}$.  
This follows from equation~(\ref{A2})
but one should also note its more detailed consequence:
$$
W_iW_j=\left( \ker \Psi_{ij}^* \circ \Theta^{(g)} \right)^{\perp}
\subseteq \left( \ker  \Theta^{(g)} \right)^{\perp}
=W_{i+j}.  
$$

For example, the even self-dual lattices of rank~$16$,
$\{ E_8 \oplus E_8, D_{16}^{+}\}$, give the
{\it problem of Witt\/}:  
find the dependence of the theta series in each genus.  
{From} results of Igusa and Kneser, cast in the above form, we have
$W_0   =W_1  =W_2  =W_3  =\<1 \>   $ and
$W_4   =V^{*} $ where $1$ is the vector of all ones.  
This is a nice way to present the linear relations.  
By a result of Igusa \cite{Ch},
$J = \vartheta^{(4)}_{6}- \vartheta^{(4)}_{7}$ gives the Schottky form in genus~$4$.  
The representation numbers for $D_4$ follow from:
$r(D_{\ell},D_4)=1152 \binom{\ell}{4}$,
$r(A_{\ell},D_4)=0$, $r(E_{8},D_4)=1152\cdot 3150$, $r(E_{7},D_4)=1152\cdot 315$,
see \cite{er}.  

\begin{thm}
\label{short1}
For $V=\C^8$, let $\Theta^{(g)}: V \to  \left[ \Gamma_g(1,2), 8 \right]$
be defined by                  
$\Theta^{(g)}(r)=\sum_{i=0}^7 r_i \vartheta_i^{(g)}$ for the
eight self-dual lattices of rank~$16$.  
For $c,\Xi,\sigma \in V^*$ given by
\begin{align*}
&\Xi=w_1(\tfrac12)=(0,2,4,8,16,32,0,0), \\
&\sigma=w_4\left(\tfrac12 D_4\right)=1152\,(140,0,630,496,3220,1820,6300,1820),  \\
&c=(0,0,0,0,0,0,1,1),
\end{align*}
the filtration $W_g= \left( \ker \Theta^{(g)} \right)^{\perp}$ is given by
$W_0   =\<1\> $, $W_1=\<1, c, \Xi\>$,  $W_2  =\<1, c, \Xi,\Xi^2\>$,
$W_3  =\<1, c, \Xi,\Xi^2,\Xi^3\>$, $W_4  =\<1, c, \Xi,\Xi^2,\Xi^3, \Xi^4, \sigma\>$
and $W_5=V^*$.  
The relation among the theta series in $g=4$ is
$\det( \vartheta^{(4)},\sigma,\Xi^4, \Xi^3, \Xi^2, \Xi , c,1)=0$.  
For the six odd lattices alone, the corresponding filtration is
$W_0   =\<1\> $, $W_1=\<1,  \Xi\>$,  $W_2  =\<1,  \Xi,\Xi^2\>$,
$W_3  =\<1, \Xi,\Xi^2,\Xi^3\>$  
and $W_4=V^*$.  
\end{thm}
{\it Proof.}
{From} the definition of $\Theta^{(0)}$, we see that $W_0=\< 1 \>$.  
{From} Table~2, we see that
$\vartheta^{(1)}=\Xi\,\Xi^{(1)}[0]+(1-c)\vartheta_0^{(1)}+c\,\vartheta_6^{(1)}$
so that $W_1=\< 1, c, \Xi \>$.  
By Theorem~\ref{T2}, the forms vanishing on the reducible locus
$\H_1 \times \H_1$ are trivial, so that
$W_2=W_1W_1=\<1, c, \Xi,\Xi^2\>$.  
By Theorem~\ref{T2}, the forms vanishing on the reducible locus
$\H_1 \times \H_2$ are trivial, so that
$W_3=W_1W_2=\<1, c, \Xi,\Xi^2,\Xi^3\>$.  
Let $r \in V$ and note $W_1W_3= \<1, c, \Xi,\Xi^2,\Xi^3,\Xi^4\>$.  
For $r \perp W_1W_3$, $\Theta^{(4)}(r)$  vanishes on the
reducible locus $\H_1 \times \H_3$ and is hence a multiple of $J^{(4)}$
by Theorem~\ref{T2}.   Thus $\Theta^{(4)}(r)= r \cdot \sigma  J^{(4)}$ by
looking at the Fourier coefficient for $\tfrac12 D_4$;
therefore  $W_4  =\<1, c, \Xi,\Xi^2,\Xi^3, \Xi^4, \sigma\>$ and the
relation in $g=4$ follows immediately.  
We have $W_5 \supseteq W_1W_4= \<1, c, \Xi,\Xi^2,\Xi^3, \Xi^4, \sigma, c\sigma, \Xi\sigma\>=V^*$.  
The corresponding result for the six odd lattices follows by restriction to the
first six coordinates.  \qed \smallskip
\begin{rem}
Hence, $\dim [\Gamma_g(1,2),8]^\vartheta $ is $3,5,6,7,8$ for
$g=1,2,3,4,5$.  
\end{rem}
It will  be important to compute
Witt images of  bases for $ [\Gamma_g(1,2),8]^\vartheta $.  For brevity, let
$$
c_0=\frac{1}{5160960}\,
\frac{\det (\sigma,\Xi^4,\Xi^3,\Xi^2,\Xi,1)}{\det (\Xi^5,\Xi^4,\Xi^3,\Xi^2,\Xi,1)}
=\frac{89 \cdot 227}{2^{19} \cdot 3\cdot 5\cdot 7^2}.  
$$
\begin{prop}
\label{short2}
Let $\{ \Xi_j \}_{j=0}^5 \subseteq \C^6$ be the dual basis to
$\{\Xi^j \}_{j=0}^5 \subseteq \C^6$.  
Write $\Theta_j^{(g)}= \Theta^{(g)} ( \Xi_j )$.  
For $g \le 4$, we have $ \Theta_g^{(g)} \in  [\Gamma_g(1,2),8]_0$.  
We have $ \Theta_0^{(g)} =\vartheta_0^{(g)} $ and
$ \Theta_5^{(4)} =c_0 J^{(4)}$.  
We have the Witt images
\begin{align*}
\Psi_{14}^{*} \Theta_5^{(5)} &=\Theta_1^{(1)} \otimes \Theta_4^{(4)} {+}
\left( 62 \Theta_1^{(1)}{+}\vartheta_0^{(1)}\right)  {\otimes} c_0 J^{(4)};  \,\,
\Psi_{23}^{*} \Theta_5^{(5)}  =\Theta_2^{(2)} \otimes \Theta_3^{(3)}, \\
\Psi_{13}^{*} \Theta_4^{(4)} &=\Theta_1^{(1)} \otimes \Theta_3^{(3)};  \quad
\Psi_{22}^{*} \Theta_4^{(4)} =\Theta_2^{(2)} \otimes \Theta_2^{(2)}, \\
\Psi_{13}^{*} \Theta_3^{(4)} &=\Theta_1^{(1)} \otimes \Theta_2^{(3)} +
\vartheta_0^{(1)} \otimes  \Theta_3^{(3)},  \\
\Psi_{13}^{*} \Theta_2^{(4)} &=\Theta_1^{(1)} \otimes \Theta_1^{(3)} +
\vartheta_0^{(1)} \otimes  \Theta_2^{(3)},  \\
\Psi_{13}^{*} \Theta_1^{(4)} &=\Theta_1^{(1)} \otimes \vartheta_0^{(3)} +
\vartheta_0^{(1)} \otimes  \Theta_1^{(3)}.  
\end{align*}
\end{prop}
{\it Proof. }  
Consider the filtration of Theorem~\ref{short1} for the six odd lattices.  
For $g \le 4$, we have
$ \Theta_g^{(g)} \in  [\Gamma_g(1,2),8]_0$ because
$\Xi_g$ is annihilated by
$\<1, \Xi, \dots, \Xi^{g-1}\>=W_{g-1}$.  
The relation $ \Theta_5^{(4)} =c_0 J^{(4)}$ follows from the $g=4$ relation in
Theorem~\ref{short1} but we may also argue directly:  
$ \Xi_5$ is annihilated by $\<1, \Xi, \dots, \Xi^{4}\>=W_1W_3$ and so
$ \Theta_5^{(4)}$ is a form vanishing on the reducible locus $\H_1 \times \H_3$,
necessarily $ \Theta_5^{(4)} =c\,J^{(4)}$ for some constant $c$ by Theorem~\ref{T2}.  
By Cramer's rule we have
$$
\Theta_5^{(4)} =
\det( \vartheta^{(4)}, \Xi^4, \Xi^3, \Xi^2, \Xi , 1)/
\det(\Xi^5, \Xi^4, \Xi^3, \Xi^2, \Xi , 1) =c\,J^{(4)}.  
$$
Evaluating at the Fourier coefficient for $\tfrac12 D_4$,  we have
$5160960 c= \det( \sigma, \Xi^4, \dots, 1)/
\det(\Xi^5, \Xi^4, \dots , 1) $ so that $c=c_0$.  
The identity $ \Theta_0^{(g)} =\vartheta_0^{(g)} $ follows from $\tau_0=0$.  

We now consider the Witt images.  
Write the map $\Theta^{(4)} \in V^{*} \otimes [\Gamma_g(1,2),8]$
in the basis $\{\Xi^j \}_{j=0}^5$ so that
$\vartheta^{(4)}= \sum_j \Theta^{(4)}(\Xi_j) \Xi^j$ or
$$
\vartheta^{(4)} = \Xi^5\, c_0J^{(4)}
+ \Xi^4\, \Theta_4^{(4)}+ \Xi^3\, \Theta_3^{(4)}+ \Xi^2\, \Theta_2^{(4)}+ \Xi\,  \Theta_1^{(4)}+ 1\, \vartheta_0^{(4)} .  
$$
That the Witt images of $\Psi_{13}^{*}$ are as stated follows from
\begin{align*}
&\Xi^4\, \Psi_{13}^{*}\Theta_4^{(4)}+\Xi^3\, \Psi_{13}^{*}\Theta_3^{(4)}+ \Xi^2\, \Psi_{13}^{*}\Theta_2^{(4)}
+ \Xi\,  \Psi_{13}^{*}\Theta_1^{(4)}+ 1\, \Psi_{13}^{*}\vartheta_0^{(4)} \\
= &\vartheta^{(1)} \otimes \vartheta^{(3)}  \\
=&\left( \Xi\,  \Theta_1^{(1)}+ 1\, \vartheta_0^{(1)} \right) \otimes
\left(  \Xi^3\, \Theta_3^{(3)}+ \Xi^2\, \Theta_2^{(3)}+ \Xi\,  \Theta_1^{(3)}+ 1\, \vartheta_0^{(3)}  \right) \\
= & \Xi^4 \left( \Theta_1^{(1)} \otimes \Theta_3^{(3)}\right)
+\Xi^3 \left( \Theta_1^{(1)} \otimes \Theta_2^{(3)} + \vartheta_0^{(1)} \otimes  \Theta_3^{(3)}\right)
+\\
\Xi^2  &\left( \Theta_1^{(1)} {\otimes} \Theta_1^{(3)} {+} \vartheta_0^{(1)} {\otimes}  \Theta_2^{(3)}\right)
{+}\Xi   \left( \Theta_1^{(1)} \otimes  \vartheta_0^{(3)} {+} \vartheta_0^{(1)} \otimes  \Theta_1^{(3)}\right)
{+} 1 \vartheta_0^{(1)} {\otimes}  \vartheta_0^{(3)}.  
\end{align*}
The others are similar but one needs to use
$\Xi^6=62 \Xi^5 -1240 \Xi^4 +9920 \Xi^3 -31744 \Xi^2 + 32768 \Xi$.  
\qed

The splitting of these forms may be used to provide finer information.  
\begin{thm}
\label{T3}  We have
$
\left[ \Gamma_4(1,2), 8 \right]_0=\C J + \C \,\Xi^{(4)}[0]  
$.  
\end{thm}
{\it Proof.}  
Take $f \in  \left[ \Gamma_4(1,2), 8 \right]_0$ and let
$\Psi_{1,3}^{*}  f= \a \,\Xi^{(1)}[0] \otimes \Xi^{(3)}[0]$.    
So $f -\a \Xi_4^{(4)}[0]$ is in $ \ker \Psi_{1,3}^{*}$ and is a multiple of $J$ by
Theorem~\ref{T2}.  \qed
\smallskip

We wish to compute the dimension of
$ \left[ \Gamma_g(1,2), 8 \right]$ for $g \le 4$.  
We know that $\dim  \left[ \Gamma_1(1,2), 8 \right]=3$,
spanned by $\vartheta_0^{(1)}$, $\Xi^{(1)}[0]$ and $\vartheta_6^{(1)}$.  
Our method for $g=3,4$ does not succeed in $g=2$, so we must make use of the result of Igusa
and  Runge that $\dim  \left[ \Gamma_2(1,2), 8 \right]=4$.  
A basis for  $ \left[ \Gamma_2(1,2), 8 \right]$
is then given by
$\vartheta_0^{(2)}$, $\Theta_1^{(2)}$,  $\Xi^{(2)}[0]$ and $\vartheta_6^{(2)}$.

\begin{thm}
\label{N2}
We have
$ \left[ \Gamma_3(1,2), 8 \right] =  \left[ \Gamma_3(1,2), 8 \right] ^{\vartheta}$
and the dimension is $5$.  
We have
$ \left[ \Gamma_4(1,2), 8 \right] =  \left[ \Gamma_4(1,2), 8 \right] ^{\vartheta}$
and the dimension is $7$.  
\end{thm}
{\it Proof.}  
We omit the proof of this theorem for $g=3$, since the fact is known
and the method that we  use is
illustrated by the $g=4$ case.  
We just observe that  a basis of
$ \left[ \Gamma_3(1,2), 8 \right] ^{\vartheta}$
is given by
$\vartheta_0^{(3)}$, $\Theta_1^{(3)}$, $\Theta_2^{(3)}$,  $\Xi^{(3)}[0]$ and $\vartheta_6^{(3)}$.  
We make use of the commutative diagram:  
\begin{align*}
&\left[ \Gamma_4(1,2), k \right]  
\quad\overset{\Psi_{13}^{*}}{\longrightarrow} \hskip1.0in
\left[ \Gamma_1(1,2), k \right] \otimes \left[ \Gamma_3(1,2), k \right] \\
&\Psi_{112}^{*}\quad \downarrow  
\hskip2.0in
\operatorname{Id} \oplus \Psi_{12}^{*}\quad \downarrow \\
&\Sym\left(  \left[ \Gamma_1(1,2), k \right]^{\otimes 2}  \right)
\otimes \left[ \Gamma_2(1,2), k \right]
\to
\left(  \left[ \Gamma_1(1,2), k \right]^{\otimes 2}  \right)
\otimes \left[ \Gamma_2(1,2), k \right]  
\end{align*}

A basis of
$ \left[ \Gamma_4(1,2), 8 \right] ^{\vartheta}$
is given by
$\vartheta_0^{(4)}$, $\Theta_1^{(4)}$, $\Theta_2^{(4)}$,  
$\Theta_3^{(4)}$,  $\Theta_4^{(4)}$,  $\Theta_5^{(4)}$  
and $\vartheta_6^{(4)}$.  
{From} Proposition~\ref{short2},
the images of $\Psi_{13}^{*}$ are
$\vartheta_0^{(1)}\otimes \vartheta_0^{(3)} $,
$\Xi^{(1)}[0] \otimes  \vartheta_0^{(3)} +\vartheta_0^{(1)} \otimes  \Theta_1^{(3)}$,
$\Xi^{(1)}[0] \otimes \Theta_1^{(3)} +\vartheta_0^{(1)} \otimes  \Theta_2^{(3)} $,  
$\Xi^{(1)}[0] \otimes \Theta_2^{(3)} +\vartheta_0^{(1)} \otimes  \Xi^{(3)}[0]$,  
$\Xi^{(1)}[0] \otimes  \Xi^{(3)}[0]$, $0$ and
$\vartheta_6^{(1)}\otimes \vartheta_6^{(3)} $.  
These linearly dependent images span a $6$~dimensional space inside
$\left[ \Gamma_1(1,2), 8 \right] \otimes \left[ \Gamma_3(1,2), 8 \right] $.  
This shows that
$\Psi_{13}^{*} \left[ \Gamma_4(1,2), 8 \right]^{\vartheta}$ is
$6$~dimensional.  

On the other hand, the general element of
the $15$~dimensional space
$\left[ \Gamma_1(1,2), 8 \right] \otimes \left[ \Gamma_3(1,2), 8 \right] $ is
given by
\begin{align*}
&\a_1 \vartheta_0^{(1)}  \otimes \vartheta_0^{(3)} +
\a_2 \Xi^{(1)}[0]                 \otimes \vartheta_0^{(3)} +
\a_3 \vartheta_6^{(1)}    \otimes \vartheta_0^{(3)} +    \\
&\b_1 \vartheta_0^{(1)}  \otimes \Theta_1^{(3)} +
\b_2 \Xi^{(1)}[0]                 \otimes \Theta_1^{(3)} +
\b_3 \vartheta_6^{(1)}    \otimes \Theta_1^{(3)}+    \\  
&\gamma_1 \vartheta_0^{(1)}  \otimes \Theta_2^{(3)} +
\gamma_2 \Xi^{(1)}[0]                 \otimes \Theta_2^{(3)} +
\gamma_3 \vartheta_6^{(1)}    \otimes \Theta_2^{(3)}+    \\  
&\delta_1 \vartheta_0^{(1)}  \otimes  \Xi^{(3)}[0] +
\delta_2 \Xi^{(1)}[0]                 \otimes  \Xi^{(3)}[0] +
\delta_3 \vartheta_6^{(1)}    \otimes  \Xi^{(3)}[0] +    \\  
&\epsilon_1 \vartheta_0^{(1)}  \otimes \vartheta_6^{(3)} +
\epsilon_2 \Xi^{(1)}[0]                 \otimes \vartheta_6^{(3)} +
\epsilon_3 \vartheta_6^{(1)}    \otimes \vartheta_6^{(3)} .
\end{align*}

By Proposition~\ref{short2},
the image of this element under
$\operatorname{Id} \oplus \Psi_{12}^{*}$   is
\begin{align*}
&\left( \a_1 \vartheta_0^{(1)}  +
\a_2 \Xi^{(1)}[0]                +
\a_3 \vartheta_6^{(1)}  \right)  \otimes  \vartheta_0^{(1)} \otimes \vartheta_0^{(2)}  +    \\
&\left( \b_1 \vartheta_0^{(1)} +
\b_2 \Xi^{(1)}[0]                +
\b_3 \vartheta_6^{(1)}  \right)  \otimes  \left( \Xi^{(1)}[0]  \otimes \vartheta_0^{(2)}
+\vartheta_0^{(1)} \otimes  \Theta_1^{(2)} \right)+    \\  
&\left( \gamma_1 \vartheta_0^{(1)} +
 \gamma_2 \Xi^{(1)}[0]                +
 \gamma_3 \vartheta_6^{(1)}  \right)  \otimes  \left( \Xi^{(1)}[0]  \otimes \Theta_1^{(2)}
+\vartheta_0^{(1)} \otimes  \Xi^{(2)}[0]  \right)+    \\  
&\left( \delta_1 \vartheta_0^{(1)}  +
\delta_2 \Xi^{(1)}[0]                +
 \delta_3 \vartheta_6^{(1)} \right)   \otimes   \Xi^{(1)}[0] \otimes  \Xi^{(2)}[0] +    \\  
&\left( \epsilon_1 \vartheta_0^{(1)}    +
\epsilon_2 \Xi^{(1)}[0]                  +
\epsilon_3 \vartheta_6^{(1)} \right)   \otimes  \vartheta_6^{(1)}  \otimes \vartheta_6^{(2)}  .
\end{align*}
If we demand that this image lie in
$\Sym\left(  \left[ \Gamma_1(1,2), 8 \right]^{\otimes 2}  \right)
\otimes \left[ \Gamma_2(1,2), 8 \right]$,
it imposes certain linear equations on the coefficients.  
Again, every term
is a tensor of basis elements.  
The free parameters are $\a_1$, $\delta_2$ and $\epsilon_3$.  
We have $\a_2=\b_1$, $\b_2=\gamma_1$ and $\gamma_2=\delta_1$.  
We have $\a_3=\b_3=\gamma_3=\delta_3=\epsilon_1=\epsilon_2=0$.  
Thus the preimage
$X=\left( \operatorname{Id} \oplus \Psi_{12}^{*} \right)\inv
\left(  \Sym\left(  \left[ \Gamma_1(1,2), 8 \right]^{\otimes 2}  \right)
\otimes \left[ \Gamma_2(1,2), 8 \right]\right)$ is $6$~dimensional inside  
$\left[ \Gamma_1(1,2), 8 \right] \otimes \left[ \Gamma_3(1,2), 8 \right] $.  
This preimage $X$ necessarily contains
$\Psi_{13}^{*} \left[ \Gamma_4(1,2), 8 \right]$.  
However, since $\Psi_{13}^{*} \left[ \Gamma_4(1,2), 8 \right]^{\vartheta}$ is
$6$~dimensional we also have
$\Psi_{13}^{*} \left[ \Gamma_4(1,2), 8 \right]^{\vartheta}=X$
$=\Psi_{13}^{*} \left[ \Gamma_4(1,2), 8 \right]$.  

{From}
$\Psi_{13}^{*} \left[ \Gamma_4(1,2), 8 \right]^{\vartheta}
=\Psi_{13}^{*} \left[ \Gamma_4(1,2), 8 \right]$
and the knowledge of the cusp forms,
we can easily deduce
$ \left[ \Gamma_4(1,2), 8 \right] =  \left[ \Gamma_4(1,2), 8 \right]^{\vartheta}$.  
For each $f \in  \left[ \Gamma_4(1,2), 8 \right] $,
there is a $g \in  \left[ \Gamma_4(1,2), 8 \right]^{\vartheta}$
with $\Psi_{13}^{*}f=\Psi_{13}^{*}g$.  
We have $\Psi_{13}^{*}(f-g)=0$ so that
$f-g\in  \left[ \Gamma_4(1,2), 8 \right]_0$.  
Thus $f= g+ \a\,  \Xi^{(4)}[0] +\b J^{(4)}\in  \left[ \Gamma_4(1,2), 8 \right]^{\vartheta}$.  
\qed

\begin{lm}
\label{M6}
Let $ f\in \left[ \Gamma_g(1,2), k\right]$ be such that
$ \Phi(f), \Phi(f|M)\in \left[ \Gamma_{g-1}, k\right]$.  
Then $\Phi\left( \Tr(f) \right)=
2^{g-1}(1+2^{g-1})\, \Phi(f) + 2^{2g-2}\, \Phi(f|M) $.  
\end{lm}
{\it Proof.}  
The trace of $f$ is given by
$
\Tr(f) = \sum_{\text{ even } \zeta} f |\gamma_{\zeta}  
$.  
According to Proposition~\ref{B4}, when
$\zeta=\begin{bmatrix}{a}\\{b}\end{bmatrix}$ and the first entry of $a$ is zero, we have
$ \Phi( f | \gamma_{\zeta})= \Phi( f | M_a \, t(S))= \Phi( f )  |  M_{a'} \, t(\pi(S))$.  
If $\Phi(f)$ is level one then $ \Phi( f | \gamma_{\zeta})=\Phi(f)$.  
There are $2^{g-1}(1+2^{g-1})$ even characteristics with the first entry of $a$  zero.

When the first entry of $a$ is one, we have
$ \Phi( f | \gamma_{\zeta})= \Phi( f | M_a \, t(S))=\Phi(f | M) | M_{c'}\, t(\pi(S))$.  
If $\Phi(f|M)$ is level one then $ \Phi( f | \gamma_{\zeta})=\Phi(f|M)$.  
There are $2^{2g-2}$ even characteristics where the first entry of $a$  is one.
Thus,  $\Phi\left( \Tr(f) \right)=
2^{g-1}(1+2^{g-1})\, \Phi(f) + 2^{2g-2}\, \Phi(f|M) $.  \qed

\begin{cor}
\label{M7}
If  $\Theta_5^{(4)}=c_0J^{(4)}$, then  
$\Tr\, \Theta_5^{(5)}= 16\cdot 17 c_0J^{(5)}$.  
\end{cor}
%
{\it Proof.}  
We know that $\Phi( \Theta_5^{(5)} ) = \Theta_5^{(4)}=c_0J^{(4)}$ is level one.  
We will show that $\Phi( \Theta_5^{(5)}|M ) = 0$.  
By Proposition~\ref{short2}, we have
$\Psi_{14}^{*} \Theta_5^{(5)} =\Theta_1^{(1)} \otimes \Theta_4^{(4)} +
( 62\,\Theta_1^{(1)}+\vartheta_0^{(1)}) \otimes c_0 J^{(4)}$.
Using that $\Theta_1^{(1)} $ is a cusp form,
we have $\Phi( \Theta_5^{(5)} |M) =\Phi( \vartheta_0^{(1)}|M ) c_0 J^{(4)}|M=0$
since $ \vartheta_0^{(1)}$ vanishes at TOC.  
By Lemma~\ref{M6}, we have
$$
\Phi( \Tr\,  \Theta_5^{(5)} )= 2^{g-1}(1+2^{g-1})\,  \Theta_5^{(4)}
= 16\cdot 17  c_0J^{(4)}.  
$$
So $ \Tr\,  \Theta_5^{(5)}\in  \left[ \Gamma_5, 8\right]$ has $\Phi$ image
$16\cdot 17 c_0J^{(4)}$.  Therefore, since \cite{PY07},  page~{216}, tells us that
the only cusp forms in $ \left[ \Gamma_5, 8\right]$ are trivial, we have
$ \Tr\,  \Theta_5^{(5)}=16\cdot 17 c_0J^{(5)}$.   \qed

\begin{df}
Let $f \in \left[ \Gamma_g(1,2), k \right]$.  
We say that $f$ is a cusp form on the Jacobian locus if
for all $\gamma \in \Gamma_g$, $\Phi(f| \gamma)$ vanishes
upon restriction to the period matrices of
compact Riemann surfaces.
\end{df}

We are ready for the  
{\it Proof of Theorem~\ref{main}. }
The map $\Theta^{(g)}: V \to [\Gamma_g(1,2),8]$ is
written as $\vartheta^{(g)}$ in the standard basis and as
$\sum_j \Theta_j^{(g)} \Xi^j$ in the $\{\Xi^j\}$ basis so that we have equation~(\ref{WE}).  
We set $$\Xi^{(g)}[0]=\Theta^{(g)}_g-
\frac{17 \cdot 89 \cdot 227}{2^{19} \cdot 3\cdot 5\cdot 7^2\cdot 33}
J^{(g)} \in [\Gamma_g(1,2),8].$$  

Since the form $ J^{(g)}$ vanishes along $\Jac_k\times \Jac_{g-k}$ when $g\leq 5$,  as an immediate consequence of Proposition~\ref{short2},  we get the splitting property for
$\Theta^{(g)}_g$ and hence for $\Xi^{(g)}[0]$ along $\Jac_k\times \Jac_{g-k}$.
Always according to  the  same Proposition, we have that the  forms $\Theta_g^{(g)}$ are cusp forms when we restrict to the Jacobian locus.
Hence, when $g\leq 4$,  their  trace is $0$ whenever we restrict to $\Jac_g$. The extra  contribution coming from $ J^{(g)}$, with $c=c_0$,  is added to get property~$(2)$ along $\H_g$.
According   to Corollary~\ref{M7}, $\Xi^{(5)}[0]$  verifies property~$(2)$ along $\Jac_5$.
Moreover,  $\Xi^{(5)}[0]$ is the unique linear combination of theta series  in
$[\Gamma_5(1,2),8]$ that satisfies both properties~$(1)$
and~$(2)$.  The  cusp forms on $\Jac_5$ from theta series are spanned by $\Xi^{(5)}[0]$ and $J^{(5)}$
and any $\Xi^{(5)}[0]+ c\,J^{(5)}$ satisfies properties~$(1)$ but,
since $J^{(5)}$ does not vanish identically on $\Jac_5$, see \cite{GSM2},
only $\Xi^{(5)}[0]$ also satisfies $(2)$.    \qed


 \section{Appendix: Theta series with harmonic polynomial coefficients}
 
A different expression for the $\Xi^{(g)}[0]$, when $ g\leq 4$,
was obtained from
theta series with harmonic polynomial coefficients. We briefly recall the results:
Let $X$ be a matrix with $m$ rows and $g$ columns.  A harmonic form of weight $\nu$ in the matrix variable  $X$ is a polynomial
$P(X) $ with the properties
\begin{equation*}
\label{H1}
\forall A\in \GL(n,\C), \quad
P(XA)=({\rm det}A)^{\nu} P(X),\notag
\end{equation*}
\begin{equation}
\label{H2}
\Delta P=\sum_{i,\, j}\frac{\partial ^2}{(\partial X_{ ij})^2}P=0.  
\notag
\end{equation}
It can be proved, cf. \cite{Fr1},  page 51 or  \cite{An}, that if $S$ is a positive definite integral unimodular quadratic form of degree
$m$ with $8$ dividing  $m$,  then the theta series
\begin{equation}
\label{Th1}
\vartheta_{S,\, P}^{(g)}(\tau) =\sum_{X\in \Z^{m,\,g} }
P(S^{1/2} X)\ e(1/2\cdot {\rm tr}(S[X]\tau))
\notag
\end{equation}
is a modular form of weight $m/2+{\nu} $ relative to $ \Gamma_g(1,2)$. It is a cusp form if ${\nu} >0$.
Moreover,  if $S$ is also even then we get a modular form relative to $\Gamma_g$.  

A simple way to construct harmonic polynomials is the  following:
let $L$  be a $m\times g$ matrix  with
$L' L=0 $ and $ L'  \overline{L}>0$,  
then, for $\nu \in \Z_{\ge 0}$,
$$P_{\nu} (X)= {\rm det} (L'X)^{\nu} $$
is a harmonic polynomial of degree $\nu$.
Here  $\overline{L}$ is the conjugate of $L$ and necessarily $m\geq 2g$.\smallskip

For $m=8$, $\nu=4$, $k=8$, $S = I_g \text{ or } E_8$ and
for $g=1,2,3,4$ we choose $L$ of the  form
$L'= \pmatrix 1_g & 0 & i1_g &0 \endpmatrix$.

\begin{prop}Let $m=8$ and $\nu=4$, then
\item{1.} The theta series
$\vartheta^{(g)} _{E_8,\, P_4}$
vanish when $g=1,2,3$ and, up to a nonzero multiplicative constant, is equal to $J$ when $g=4$. \smallskip
 
\item{2.} The theta series
$\vartheta^{(g)}_{1_8,\, P_4}$
do not vanish when $g=1,2,3,4$.
\end{prop}
{\it Proof}. In \cite{Ho} and \cite{W1} the non-vanishing
of $\vartheta^{(4)}_{E_8,\, P_4}$ has been proved.  
The vanishing of the other cases is a consequence of the general  
fact that there are no level one cusp forms of weight $8$ when $g\leq 3$.
About the second statement, we observe that the nonvanishing is the consequence of a simple computation, since in all
these cases, we have $a(\tfrac12 1_g)\neq 0$ for the Fourier coefficients of $\tfrac12 1_g$.
In fact
$$a(\tfrac12 1_g)=\sum_{X\in \Z^{m,\,g} : \,\,\,X'X=1_g }  {\rm det} (L'X)^{4} .$$

For such $X$ it is easy to check that ${\rm det} (L'X)$ is 0 or a fourth root of unity. Since there exist $X$
such that the previous determinant is not  zero, we get $a(\tfrac12 1_g)\neq 0$. \qed

\end{document}